\documentclass[a4paper,12pt]{amsart}
\usepackage{amssymb}
\usepackage{ifthen}
\usepackage{graphicx}
\usepackage{float}
\usepackage{caption}
\usepackage{subcaption}
\usepackage{cite}
\usepackage{amsfonts}
\usepackage{amscd}
\usepackage{amsxtra}
\usepackage{color}
\usepackage{enumerate}

\newtheorem{thm}{Theorem}[section]
\newtheorem{lem}[thm]{Lemma}
\newtheorem{cor}[thm]{Corollary}

\newtheorem{cl}{Claim}[section]
\newtheorem{ca}{Case}[section]
\newtheorem{sca}{Subcase}[section]
\newtheorem{scl}[section]{Subclaim}
\newtheorem{conj}[equation]{Conjecture}

\theoremstyle{definition}
\newtheorem{defn}[thm]{Definition}
\newtheorem{prob}[equation]{Problem}
\newtheorem{op}[thm]{Open Problem}
\newtheorem{ques}[thm]{Question}
\newtheorem{rem}[thm]{Remark}
\newtheorem{exam}[equation]{Example}
\newenvironment{pf}[1][]{%
 \vskip 3mm
 \noindent
 \ifthenelse{\equal{#1}{}}%
  {{\slshape Proof. }}%
  {{\slshape #1.} }%
 }%
{\qed\bigskip}

\newcounter{alphabet}

\newcommand{\id}{{\operatorname{id}}}

\newcommand{\diam}{{\operatorname{diam}}}

\def\be{\begin{equation}}
\def\ee{\end{equation}}

\newcommand{\ben}{\begin{enumerate}}
\newcommand{\een}{\end{enumerate}}

\newcommand{\blem}{\begin{lem}}
\newcommand{\elem}{\end{lem}}
\newcommand{\bthm}{\begin{thm}}
\newcommand{\ethm}{\end{thm}}
\newcommand{\bcor}{\begin{cor}}
\newcommand{\ecor}{\end{cor}}
\newcommand{\beg}{\begin{exam}}
\newcommand{\eeg}{\end{exam}}
\newcommand{\begs}{\begin{examples}}
\newcommand{\eegs}{\end{examples}}
\newcommand{\bdefe}{\begin{defn}}
\newcommand{\edefe}{\end{defn}}
\newcommand{\bprob}{\begin{prob}}
\newcommand{\eprob}{\end{prob}}
\newcommand{\bques}{\begin{ques}}
\newcommand{\eques}{\end{ques}}
\newcommand{\bei}{\begin{itemize}}
\newcommand{\eei}{\end{itemize}}
\newcommand{\bcon}{\begin{conj}}
\newcommand{\econ}{\end{conj}}
\newcommand{\bop}{\begin{op}}
\newcommand{\eop}{\end{op}}

\newcommand{\bca}{\begin{ca}}
\newcommand{\eca}{\end{ca}}
\newcommand{\bsca}{\begin{sca}}
\newcommand{\esca}{\end{sca}}

\newcommand{\bcl}{\begin{cl}}
\newcommand{\ecl}{\end{cl}}

\newcommand{\bscl}{\begin{scl}}
\newcommand{\escl}{\end{scl}}

\newcommand{\bcons}{\begin{conjs}}
\newcommand{\econs}{\end{conjs}}
\newcommand{\bprop}{\begin{propo}}
\newcommand{\eprop}{\end{propo}}
\newcommand{\br}{\begin{rem}}
\newcommand{\er}{\end{rem}}
\newcommand{\brs}{\begin{rems}}
\newcommand{\ers}{\end{rems}}
\newcommand{\bo}{\begin{obser}}
\newcommand{\eo}{\end{obser}}
\newcommand{\bos}{\begin{obsers}}
\newcommand{\eos}{\end{obsers}}
\newcommand{\bpf}{\begin{pf}}
\newcommand{\epf}{\end{pf}}
\newcommand{\ba}{\begin{array}}
\newcommand{\ea}{\end{array}}
\newcommand{\beq}{\begin{eqnarray}}
\newcommand{\beqq}{\begin{eqnarray*}}
\newcommand{\eeq}{\end{eqnarray}}
\newcommand{\eeqq}{\end{eqnarray*}}

\newcounter{minutes}
\divide\time by 60
\newcounter{hours}
\multiply\time by 60 \addtocounter{minutes}{-\time}

\begin{document}

\bibliographystyle{amsplain}
\title{A note on $\partial$-bilipschitz mappings in quasiconvex metric spaces}

\thanks{
File:~\jobname .tex,
          printed: \number\day-\number\month-\number\year,
          \thehours.\ifnum\theminutes<10{0}\fi\theminutes}

\author{Tiantian Guan}
\address{Tiantian Guan, School of Mathematics and Big Data, Foshan University,  Foshan, Guangdong 528000, People's Republic
of China} \email{ttguan93@163.com}

\author{Saminathan Ponnusamy}
\address{Saminathan Ponnusamy, Department of Mathematics, Indian Institute of Technology Madras, Chennai 600036,
India}
\address{Department of Mathematics, Petrozavodsk State University, ul., Lenina 33, 185910 Petrozavodsk,
Russia}
\email{samy@iitm.ac.in}

\author{Qingshan Zhou $^{~\mathbf{*}}$
}
\address{Qingshan Zhou, School of Mathematics and Big Data, Foshan University,  Foshan, Guangdong 528000, People's Republic of China}
\email{qszhou1989@163.com; q476308142@qq.com}

%

\date{}

\keywords{$\partial$-biLipschitz mapping, quasihyperbolic metric, quasiconvex metric space, uniform space, quasim\"obius mapping
\\
${}^{\mathbf{*}}$ Corresponding author
}

\begin{abstract}
This paper focuses on properties of $\partial$-biLipschitz mappings which were recently introduced by Bulter. We establish several characterizations for the class of $\partial$-biLipschitz mappings between domains in quasiconvex metric spaces. As an application, we show that a locally quasisymmetric equivalence between uniform metric spaces is quasim\"obius, quantitatively.
\end{abstract}

\thanks{Qingshan Zhou was supported by NNSF of China (No. 11901090), the Department of Education of Guangdong Province, China (No. 2021KTSCX116), Guangdong-Hong Kong-Macao Intelligent Micro-Nano Optoelectronic Technology Joint Laboratory (Project No. 2020B1212030010), and Guangdong Basic and Applied Basic Research Foundation (Grant no. 2021A1515012289).}

\maketitle{}
\pagestyle{myheadings}
\markboth{T. Guan, S. Ponnusamy and Q. Zhou}{A note on $\partial$-biLipschitz mappings in quasiconvex metric spaces}

\section{Introduction and main results}\label{sec-1}
In order to explain the connections between the uniformizations of Gromov hyperbolic spaces established by Bonk et al. \cite{BHK} and  Bulter \cite{Bu1} (see also \cite{ZP}), Bulter recently introduced the following class of $\partial$-biLipschitz mappings in \cite{Bu}. 

\begin{defn}
Let $L\geq 1$ and $0<\lambda<1$. A homeomorphism $f:\, \Omega\to \Omega'$ between incomplete, locally complete and connected metric spaces is called $\partial$-{\it Lipschitz} with data $(L, \lambda)$ if for each $x\in \Omega$ and for all $y, z\in B(x, \lambda d_{\Omega}(x))$,
$$\frac{d'(f(y), f(z))}{d_{\Omega'}(f(x))}\leq L\frac{d(y, z)}{d_{\Omega}(x)}.
$$
Moreover, $f$ is called $\partial$-{\it biLipschitz} with data $(L, \lambda)$ if  both $f$ and its inverse $f^{-1}$ are $\partial$-Lipschitz with data $(L, \lambda)$. Here $d_{\Omega}(x)$ denotes the distance between $x$ and the metric boundary $\partial \Omega$ for all $x\in \Omega$, and
$$B(x,\lambda d_{\Omega}(x))=\{z\in X:\; d(z,x)<\lambda d_{\Omega}(x)\}.
$$
The local completeness ensures that $d_{\Omega}(x)$ does not vanish.
\end{defn}

Several important conformal transformations are $\partial$-biLipschitz mappings (cf. \cite{Bu}). 
Examples of such mappings include Riemann mappings from the unit disc onto simply connected plane domains, sphericalization and inversion operations in metric spaces setting suggested independently by Bonk and Kleiner \cite{BK02} and Buckley et al. \cite{BHX}, and conformal deformations introduced by Bonk et al. \cite{BHK} and Bulter \cite{Bu1}. More generally, Bulter proved that a homeomorphism between uniform spaces is $\partial$-biLipschitz if and only if it is biLipschitz with respect to the quasihyperbolic metric, see \cite[Theorem 1.5]{Bu}.

This work is a continuation of \cite{Bu}. We investigate new characterizations of $\partial$-biLipschitz mappings between domains (open and connected nonempty subsets) in quasiconvex metric spaces.  For more investigations concerning properties of quasiconformal mappings in this setting, one may  refer  \cite{HL,HLL,HRWZ,ZLH}.

Obviously, the uniformity property of metric spaces implies the quasiconvexity. Note that a metric space $(X,d)$ is called {\it quasiconvex} if there is a constant $c\geq 1$ such that for each $x, y\in X$, there is a curve $\gamma$ in $X$ joining $x$ and $y$ whose length $\ell(\gamma)$ satisfies the condition $\ell(\gamma)\leq cd(x, y).$ It is worth mentioning that the class of quasiconvex spaces is very large which includes convex domains, domains $G\subset \mathbb{R}^n$ endowed with its length metric associated to the Euclidean metric, uniform spaces, Ahlfors' regular spaces admitting Loewner conditions introduced by Heinonen and Koskela \cite{HeK}, and more generally doubling metric measure spaces which satisfy Poincar\'e inequalities \cite{LS}.

First, we establish the following characterizations of $\partial$-Lipschitz mappings. The terminology appeared in Theorem \ref{thm-1} and in the rest of this section will be explained in Section \ref{sec-2}.

\begin{thm}\label{thm-1}
Suppose that $(X, d)$ and $(X', d')$ are $c$-quasiconvex and complete metric spaces, that $G\subsetneq X$ and $G'\subsetneq X'$ are domains, and that $f:G\to G'$ is a homeomorphism. Then the following statements are equivalent:
\begin{enumerate}[{\rm\;\;\;\;\;(1)}]
\item\label{en1-1} $f$ is $\partial$-Lipschitz with data $(L, \lambda)$;
\item\label{en1-2} $f$ is $(\theta, t_0)$-relative with $\theta(t)=c_1t$;
\item\label{en1-3} $f$ is $\varphi$-semisolid with $\varphi(t)=c_2t$.
\end{enumerate}
Note that the parameters depend only on each other and $c$.
\end{thm}

\br  The class of relative mappings between domains in $\mathbb{R}^n$  was introduced by Gehring in his study of quasiconformal mappings. Originally, solid mappings were investigated by Gehring and Osgood \cite{GO} in proving that quasiconformal mappings in Euclidean spaces are solid mappings and the converse is not true, see also \cite{TV2}. In \cite{Vai-5}, V\"ais\"al\"a established the theory of freely quasiconformal mappings in Banach spaces by using the quasihyperbolic metric and solid mappings. Recently, Huang et al.  \cite{HRWZ} studied the relationships between solid mappings and (weakly) quasisymmetric mappings in quasiconvex and complete metric spaces.
\er

Using Theorem \ref{thm-1}, we obtain many characterizations for the class of $\partial$-biLipschitz mappings between domains in quasiconvex metric spaces. The following theorem is  an improvement of \cite[Theorem 1.5]{Bu}, where the domains are assumed to be uniform.

\begin{thm}\label{thm-2}
Suppose that $(X, d)$ and $(X', d')$ are $c$-quasiconvex and complete metric spaces, $G\subsetneq X$ and $G'\subsetneq X'$ are domains, and that $f:G\to G'$ is a homeomorphism. Then the following statements are equivalent:
\begin{enumerate}[{\rm \;\;\;\;\;(1)}]
\item\label{en2-1} $f$ is $\partial$-biLipschitz with data $(L, \lambda)$;
\item\label{en2-2} $f$ is $M$-biLipschitz with respect to the quasihyperbolic metrics;
\item\label{en2-3} $f$ and $f^{-1}$ are $(\theta, t_0)$-relative with $\theta(t)=c_1t$;
\item\label{en2-4} $f$ and $f^{-1}$ are $(L_1, \lambda_1)$-locally biLipschitz;
\item\label{en2-5} $f$ and $f^{-1}$ are $(\eta, q)$-locally quasisymmetric with $\eta(t)=c_2t$.
\end{enumerate}
Note that the parameters depend only on each other and $c$.
\end{thm}

\br The class of locally biLipschitz mappings mentioned in Theorem \ref{thm-2}(\ref{en2-4}) was introduced by Bonk et al. \cite{BHK}, where they demonstrated that there is a one-to-one correspondence between Gromov hyperbolic space category with biLipschitz mappings as morphisms and uniform space category with quasisimilarity mappings as morphisms. For more investigations about these mappings, see \cite{Z,ZLH}.
\er

As mentioned before, Bulter used $\partial$-biLipschitz mappings to investigate the relationships between lots of conformal deformations on uniform metric spaces, such as  the uniformization  processes of Gromov hyperbolic spaces in \cite{BHK,Bu1,ZP} and the sphericalization and inversion deformations in \cite{BK02,BHX}. In particular, he obtained the following:

\vspace{8pt}
\noindent
{\bf Theorem A.} $($\cite[Theorem 1.4]{Bu}$)$ 
{\it Suppose that $(\Omega, d)$ and $(\Omega', d')$ are locally compact and $A$-uniform metric spaces. If $f:\, \Omega\to \Omega'$ is $\partial$-biLipschitz with data $(L, \lambda)$, then $f$ is $\theta$-quasim\"obius with $\theta$ depending only on $A$, $\lambda$, and $L$.
}
\vspace{8pt}

The class of uniform Euclidean domains  was independently introduced by Jones in \cite{Jo80} and by Martio and Sarvas in \cite{MS}. It is known that this notation plays an important role in geometric function theory and analysis in metric spaces, see \cite{BHK,BHX,Bu1,Vai-4,ZP} and the references therein. In \cite{Vai-5}, V\"ais\"al\"a proved that a quasiconformal mapping between uniform Euclidean domains is quasim\"obius. A well-known fact is that quasiconformal mappings and their inverses are both locally quasisymmetric.

It is natural to consider whether Theorem~A 
holds or not for the class of locally quasisymmetric equivalence mappings. We say that a homeomorphism between two incomplete and locally complete metric spaces is locally quasisymmetric equivalence if  the mapping and its inverse are both locally quasisymmetric. As an application of Theorem \ref{thm-2}, we provide an affirmative answer to this question by showing the following.

\begin{thm}\label{thm-3} Let $A\geq 1$, $0<q< 1$ and $\eta:[0, \infty)\to [0, \infty)$ a homeomorphism. Suppose that $(\Omega, d)$ and $(\Omega', d')$ are locally compact and $A$-uniform metric spaces. If a homeomorphism $f:\Omega\to \Omega'$ and $f^{-1}$ are both $(\eta, q)$-locally quasisymmetric, then $f$ is $\theta$-quasim\"obius, where $\theta$ depends only on $A$, $q$, and $\eta$.
\end{thm}

\br
$(a)$ Theorem \ref{thm-3} is a generalization of Theorem~A 
and the result of V\"ais\"al\"a in \cite{Vai-5}. Indeed, one observes from Theorem \ref{thm-2} that $\partial$-biLipschitz mappings are $(\eta, q)$-locally quasisymmetric mappings with a special control function $\eta$.

$(b)$ Our method of proof of  Theorem \ref{thm-3} is also different from that of \cite{Bu,Vai-5}. In \cite{Vai-5}, the tools of conformal modulus and $n$-Lebesgue measures are needed to show that uniform domains satisfy a quasiextremal distance condition in the sense of Gehring and Martio. In order to establish Theorem~A, 
Bulter made use of the unbounded uniformizations of Gromov hyperbolic spaces associated to a class of conformal densities via Busemann functions, see \cite{Bu1,ZP}. In particular, he proved that the identity mapping between any two deformed uniform spaces is quasim\"obius.

In this paper we do not use the results established in \cite{Bu1}, but apply a result of quasisymmetric mappings from local to global given in \cite[Theorem 1.4]{ZP} and also the metric spaces sphericalization transformations introduced in \cite{BHX}.
\er

The rest of this paper is organized as follows. In Section \ref{sec-2}, we recall necessary definitions and preliminary results. The proofs of Theorems \ref{thm-1} and \ref{thm-2} are given in Section \ref{sec-3}. Section \ref{sec-4} is devoted to the proof of Theorem \ref{thm-3}.

\section{Preliminaries and notations}\label{sec-2}
Let $(X,d)$ be a metric space. For a bounded set $E\subset X$, the diameter of $E$ is denoted by $\diam(E)$. The open (resp. closed) metric ball with center $x\in X$ and radius $r>0$ is denoted by
$$B(x,r)=\{z\in X:\; d(z,x)<r\}\;\;(\mbox{resp.}\;\; \overline{B}(x,r)=\{z\in X:\; d(z,x)\leq r\}),
$$
and the metric sphere by $S(x,r)=\{z\in X:\; d(z,x)=r\}.$
A geodesic arc $\gamma$ joining $x$ and $y$ in $X$ is a continuous map $\gamma:\, I=[0, l]\to X$ from an interval $I$ to $X$ such that $\gamma(0)=x$ and $\gamma(l)=y$ and $d(\gamma(t), \gamma(t'))=|t-t'|$ for all $t, t'\in I$. If $I=[0,\infty)$, then $\gamma$ is called a geodesic ray. We say that $X$ is geodesic if every pair of points can be joined by a geodesic arc. The space $X$ is called rectifiably connected if each pair of points in $X$ can be connected by a curve $\gamma$ with $\ell(\gamma)<\infty$, where $\ell(\gamma)$  is the length of $\gamma$.

\subsection{Quasihyperbolic metric and uniform space}
Let $(\Omega, d)$ be a rectifiably connected, locally complete and incomplete metric space, where the identity mapping $\id: (\Omega, d)\to (\Omega, \ell)$ is a homeomorphism and $\ell$ is the length metric of $\Omega$ with respect to $d$. Here $\Omega$ is incomplete if its metric boundary $\partial \Omega=\bar{\Omega}\setminus \Omega$ is non-empty, where $\bar{\Omega}$ designates the metric completion of $\Omega$.
Note that every domain $G \subsetneq X$ in a quasiconvex complete space equipped with the metric $d$ is such a  space.

For all $x, y\in \Omega$, the number
$$k_{\Omega}(x, y)=\inf_{\gamma}\int_{\gamma}\frac{ds}{d_{\Omega}(z)}$$
is called {\it the quasihyperbolic distance} between $x$ and $y$ in $\Omega$, where the infimum is taken over all rectifiable curves in $\Omega$ joining $x$ and $y$.

\begin{defn}
Let $A\geq 1$. The space $\Omega$ is called $A$-{\it uniform} if for each $x, y\in \Omega$, there is a curve $\gamma$ joining $x$ and $y$ such that
\begin{enumerate}
\item $\ell(\gamma)\leq A d(x, y)$,
\item $\min\{\ell(\gamma[x,z]), \ell(\gamma[z,y])\}\leq A d_{\Omega}(z)$,
\end{enumerate}
where $\gamma[x,z]$  is the subcurve of $\gamma$ between $x$ and $z$. The curve $\gamma$ is called an $A$-uniform curve.
\end{defn}

The following facts about the quasihyperbolic metric are needed for our purpose.

\vspace{8pt}
\noindent
{\bf Lemma B.} $($\cite[Lemma 3.8]{HRWZ}$)$
{\it
Let  $(X, d)$ be a $c$-quasiconvex and complete metric space and let $G\subsetneq X$ be a domain.
\begin{enumerate}[{\rm \;\;\;\;\;(1)}]
\item\label{HRWZ-1} For all $x, y\in G$,
$d(x, y)\leq (e^{k_G(x, y)}-1)d_G(x);$
\item\label{HRWZ-2} Suppose that $x, y\in G$ and either $d(x, y)\leq \frac{1}{3c}d_G(x)$ or $k_G(x, y)\leq 1$. Then
    $$\frac{1}{2}\frac{d(x, y)}{d_G(x)}\leq k_G(x, y)\leq 3c\frac{d(x, y)}{d_G(x)}.$$
\end{enumerate}
}
\vspace{8pt}

Given a Euclidean domain $G$ with $x\in G$, we know that the ball $B(x, r)$ lies in $G$ for $0<r<d_G(x)$. However, this result is not true in the setting of metric spaces, see \cite[Example 3.1]{HRWZ}. In \cite{HLL},  Huang et al.  established the following:

\vspace{8pt}
\noindent
{\bf Lemma C.} $($\cite[Lemma 3.4]{HRWZ}$)$
{\it
Suppose that $X$ is a $c$-quasiconvex and complete metric space and $G\subsetneq X$ is a domain. Let $x\in G$ and $0<r\leq \frac{2}{2+c}d_G(x)$. Then $B(x, r)\subset G$.
}
\vspace{8pt}

\subsection{Mappings}

\begin{defn}
Let $L\geq 1$, $C\geq 0$ and let $\eta:[0, \infty)\to [0, \infty)$ be a homeomorphism. A map $f:\, X\to X'$ between metric spaces is called
\begin{enumerate}
  \item $(L,C)$-{\it quasi-isometric} if for all $x$, $y\in X$,
$$
L^{-1}d(x, y)-C\leq d'(f(x), f(y))\leq L d(x, y)+C.
$$
If $C=0$, then $f$ is called $L$-{\it biLipschitz}.
  \item $\eta$-{\it quasisymmetric} if $f$ is a homeomorphism and for every distinct points $x, a, b\in X$, we have
$$\frac{d'(f(x), f(a))}{d'(f(x), f(b))}\leq \eta\left(\frac{d(x, a)}{d(x, b)}\right).$$
  \item $\eta$-quasim\"obius if $f$ is a homeomorphism and for any distinct points $x, y, z, w\in X$, we have
$$\frac{d'(f(x), f(z))d'(f(y), f(w))}{d'(f(x), f(y))d'(f(z), f(w))}\leq \eta\left(\frac{d(x, z)d(y, w)}{d(x, y)d(z, w)}\right).$$
\end{enumerate}
\end{defn}

\begin{defn}
Let $L\geq 1$, $0<q< 1$ and let $\eta:[0, \infty)\to [0, \infty)$ be a homeomorphism. A homeomorphism $f:\, \Omega\to \Omega'$ between incomplete, locally complete and rectifiably connected metric spaces is called
\begin{enumerate}
  \item $(\eta, q)$-{\it relative} if for each $x\in \Omega$ and for all $y\in B(x, q d_{\Omega}(x))$, we have
$$\frac{d'(f(x), f(y))}{d_{\Omega'}(f(x))}\leq \eta\left(\frac{d(x, y)}{d_{\Omega}(x)}\right).$$
  \item $\eta$-{\it semisolid} if for all $x,y\in \Omega$,
$k_{\Omega'}(f(x), f(y))\leq \eta(k_\Omega(x, y)).$
  \item $(L, q)$-{\it locally biLipschitz} if for each $x\in \Omega$, there is a constant $C_x$ (depending on the point $x$) such that for all $y, z\in B(x, q d_{\Omega}(x))$, we have
$$L^{-1}C_x d(y, z)\leq d'(f(y), f(z))\leq LC_x d(y, z).$$
  \item $(\eta, q)$-{\it locally quasisymmetric} if for every $x\in \Omega$, the restriction $f|_{B(x, qd_{\Omega}(x))}$ of $f$ on $B(x, qd_{\Omega}(x))$ is $\eta$-quasisymmetric.
\end{enumerate}
\end{defn}

Finally, we record an auxiliary result concerning the quasisymmetry from local to global between uniform spaces in a quantitative way.

\vspace{8pt}
\noindent
{\bf Lemma D.} $($\cite[Theorem 1.4]{ZP}$)$
{\it Let $A\geq 1$, $0<q<1$, and let $\eta:\,[0, \infty)\to [0, \infty)$ be a homeomorphism. Suppose that $f:\Omega\to \Omega'$ is a homeomorphism between bounded, locally compact, and $A$-uniform spaces. If $f$ and $f^{-1}$ are both $(\eta, q)$-locally quasisymmetric, then $f$ is $\eta_0$-quasisymmetric if and only if there is a number $C_0 \geq 1$ and a point $w\in \Omega$ such that
${\rm diam}(\Omega)\leq C_0d_{\Omega}(w)$ and ${\rm diam}(\Omega')\leq C_0d_{\Omega'}'(f(w)).$  The parameters $\eta_0$ and $C_0$ depend on each other and also on $\eta$, $q$, and $A$.
}

\section{Proofs of Theorems \ref{thm-1} and \ref{thm-2}}\label{sec-3}

In this section, we assume that $(X, d)$ and $(X', d')$ are $c$-quasiconvex and complete metric spaces, that $G\subsetneq X$ and $G'\subsetneq X'$ are domains, and that $f:G\to G'$ is a homeomorphism.

\subsection{Proof of Theorem \ref{thm-1}}
We divide the proof of Theorem \ref{thm-1} into several lemmas.
\begin{lem}\label{lem3-1}
Assume that $L\geq 1,$ $0<\lambda<1$, and $f$ is $\partial$-Lipschitz with data $(L, \lambda).$ Then $f$ is $(\theta, t_0)$-relative with $\theta(t)=Lt$ and $t_0=\lambda$.
\end{lem}
\bpf
Let $x\in G$ and $y\in B(x, \lambda d_G(x))$. Since $f$ is $\partial$-Lipschitz, we have
$$\frac{d'(f(x), f(y))}{d_{G'}(f(x))}\leq L\frac{d(x, y)}{d_G(x)},$$
which implies that $f$ is $(\theta, t_0)$-relative with $\theta(t)=Lt$ and $t_0=\lambda$.
\epf
\begin{lem}\label{lem3-2}
Assume that  $0<t_0\leq 1,$ $c_1\geq 1$, and $f$ is $(\theta, t_0)$-relative with $\theta(t)=c_1t.$ Then $f$ is $\varphi$-semisolid with $\varphi(t)=c_2t$, where $c_2$ depends only on $c_1$, $t_0$, and $c$.
\end{lem}
\bpf Put
$$t_1=\min\left\{\log\left(\frac{t_0}{2}+1\right), \log\left(1+\frac{1}{3c_1c}\right)\right\}\,\,\,\mbox{and}\,\,\,\,\,\,c_2= \frac{24cc_1}{t_1}.$$
We first establish the following:
\bcl\label{z-1} If  $x, y\in G$ with $k_G(x, y)\leq t_1$, then $k_{G'}(f(x), f(y))\leq 6cc_1k_G(x, y)$. \ecl
Fix $x, y\in G$ with $k_G(x, y)\leq t_1$. It follows from Lemma~B(\ref{HRWZ-1}) 
that
$$d(x, y)\leq (e^{t_1}-1)d_G(x)\leq \min\left\{\frac{t_0}{2}, \frac{1}{3cc_1}\right\}d_G(x).
$$
Since $f$ is $(\theta, t_0)$-relative with $\theta(t)=c_1t$, we have
\beq\label{eq3-1}
\frac{d'(f(x), f(y))}{d_{G'}(f(x))}\leq \theta\Big(\frac{d(x, y)}{d_G(x)}\Big)=c_1\frac{d(x, y)}{d_G(x)}\leq \frac{1}{3c}.
\eeq
Now we obtain from \eqref{eq3-1} and Lemma~B(\ref{HRWZ-2}) 
that
\beq\label{eq3-2}
k_{G'}(f(x), f(y))\leq 3c\frac{d'(f(x), f(y))}{d_{G'}(f(x))}\leq 1.
\eeq
Again using Lemma~B(\ref{HRWZ-2}), 
we see from \eqref{eq3-1} and \eqref{eq3-2} that
$$k_{G'}(f(x), f(y))\leq 3c\frac{d'(f(x), f(y))}{d_{G'}(f(x))}\leq 3cc_1\frac{d(x, y)}{d_G(x)}\leq 6cc_1k_G(x, y),$$
and therefore, Claim \ref{z-1} is proved.

Next, we show that
$k_{G'}(f(x), f(y))\leq c_2k_G(x, y)$
holds for all $x,y\in G$.

Fix $x,y\in G$. If  $k_G(x, y)\leq t_1$, then the required assertion follows from Claim \ref{z-1}. So it suffices to consider the case that $k_G(x, y)> t_1$. By \cite[Lemma 3.9]{HRWZ}, the space $(G, k_G)$ is $2$-quasiconvex. Therefore, there exists a curve $\gamma$ in $G$ joining $x$ and $y$ such that
\beq\label{eq3-3}
\ell_{k_G}(\gamma)\leq 2k_G(x, y),
\eeq
where $\ell_{k_G}(\gamma)$ denotes the length of $\gamma$ in the quasihyperbolic metric.
Take a sequence of points $x=x_0$, $x_1$, $x_2$, $\ldots$, $x_{n-1}$ in $\gamma$ satisfying
$k_G(x_i, x_{i-1})=t_1$ and $k_G(x_{n-1}, y)\leq t_1,$
for $i=1, 2, \ldots, n-1$, and denote $x_n=y$. Then we know from Claim \ref{z-1} that $n\geq 2$ and $k_{G'}(f(x_i), f(x_{i-1}))\leq  6cc_1t_1\leq 6cc_1$
for $i=1, 2, \ldots, n$. Thus, one finds that
\beq\label{eq3-3a}
k_{G'}(f(x), f(y))\leq \sum\limits_{i=1}^{n}k_{G'}(f(x_i), f(x_{i-1}))\leq 6cc_1 n\leq 12cc_1(n-1).
\eeq
It follows from the inequality \eqref{eq3-3} that $(n-1) t_1\leq \ell_{k_G}(\gamma)\leq 2k_G(x, y),$
which, together with \eqref{eq3-3a}, implies $k_{G'}(f(x), f(y))\leq 24cc_1t_1^{-1}k_G(x, y).$ The proof is complete.
\epf

\begin{lem}\label{lem3-3}
If $f$ is $\varphi$-semisolid with $\varphi(t)=c_2t$ for some  $c_2\geq 1$, then $f$ is $\partial$-Lipschitz with data $(L, \lambda)$, where $L$ and $\lambda$ depend only on $c_1$ and $c$.
\end{lem}
\bpf
Let $\lambda=1/(36c^2c_2)$ and $L=24cc_2$. For  all $x\in G$ and $y, z\in B(x, \lambda d_G(x))$, we have by the triangle inequality
\beqq\label{eq3-4}
d(y, z)\leq \frac{1}{18c^2c_2}d_G(x),
\eeqq
and
\beq\label{eq3-5}
d_G(y)\geq d_G(x)-d(x, y)\geq \frac{1}{2}d_G(x).
\eeq
Therefore, we have
$d(y, z)\leq d_G(y)/(9c^2c_2)$.  Using Lemma~B(\ref{HRWZ-2}), 
we  obtain
\beq\label{eq3-7}
k_G(y, z)\leq 3c\frac{d(y, z)}{d_G(y)}\leq \frac{1}{3cc_2}.
\eeq
Because $y\in B(x, \lambda d_G(x))$, Lemma~B(\ref{HRWZ-2}) 
leads to
\beqq\label{eq3-8}
k_G(x, y)\leq 3c\frac{d(x, y)}{d_G(x)}\leq \frac{1}{12cc_2}.
\eeqq
Since $f$ is $\varphi$-semisolid,  we see from \eqref{eq3-7} that
$k_{G'}(f(x), f(y))\leq c_2k_G(x, y)\leq 1/(12c)$
and
$k_{G'}(f(y), f(z))\leq c_2k_G(y, z)\leq 1/(3c).$
Together with Lemma~B(\ref{HRWZ-2}), 
these two inequalities imply that
\beq\label{eq3-9}
\frac{d'(f(y), f(z))}{d_{G'}(f(y))}\leq 2k_{G'}(f(y), f(z)),
\eeq
and
\beqq
\frac{d'(f(x), f(y))}{d_{G'}(f(x))}\leq 2k_{G'}(f(x), f(y))\leq \frac{1}{6c}.
\eeqq
Therefore, we have
$$
d_{G'}(f(y))\leq d_{G'}(f(x))+d'(f(x), f(y))\leq (1+1/(6c))d_{G'}(f(x))\leq 2d_{G'}(f(x)).
$$
Then it follows from \eqref{eq3-5}, \eqref{eq3-7}, \eqref{eq3-9} and the semisolidity of $f$ that
\beqq
\frac{d'(f(y), f(z))}{d_{G'}(f(x))}&\leq& \frac{2d'(f(y), f(z))}{d_{G'}(f(y))}\leq 4k_{G'}(f(y), f(z))   \\
&\leq&4c_2k_G(y, z)\leq 12cc_2\frac{d(y, z)}{d_G(y)} \leq 24cc_2\frac{d(y, z)}{d_G(x)}.
\eeqq
The proof of the lemma is complete.
\epf

\subsection*{Proof of Theorem \ref{thm-1}}
The implications (\ref{en1-1})$\Rightarrow$ (\ref{en1-2})$\Rightarrow$ (\ref{en1-3})$\Rightarrow$ (\ref{en1-1}) in Theorem \ref{thm-1} follow from Lemmas \ref{lem3-1}, \ref{lem3-2}, and \ref{lem3-3}, respectively.
\qed

\subsection{Proof of Theorem \ref{thm-2}}
Before the proof, we establish some results.

\begin{lem}\label{lem3-4}
Let $0<t_0\leq 1$ and $c_1\geq 1$. If both $f$ and $f^{-1}$ are $(\theta, t_0)$-relative with $\theta(t)=c_1t$, then both $f$ and the inverse $f^{-1}$ are $(L_1, \vartheta_1)$-locally biLipschitz, where $L_1$ and $\vartheta_1$ depend only on $t_0$ and $c_1$.
\end{lem}
\bpf Let $\vartheta_1=t_0/(8c_1)$ and $L_1=4c_1$. We only need to show that $f$ is $(L_1, \vartheta_1)$-locally Lipschitz, because the other direction for $f^{-1}$ follows from a similar argument. To this end, for all $x\in G$ and $y, z\in B(x, \vartheta_1d_G(x))$, we have
\beq\label{eq4-1}
\frac{1}{2}d_G(x)\leq d_G(y)\leq 2d_G(x),
\eeq
and
$$d(y, z)\leq 2\vartheta_1d_G(x)\leq \frac{2\vartheta_1}{1-\vartheta_1}d_G(y)<t_0d_G(y).$$
Since $f$ is $(\theta, t_0)$-relative, we obtain
\beq\label{eq4-2a}
\frac{d'(f(y), f(z))}{d_{G'}(f(y))}\leq c_1\frac{d(y, z)}{d_G(y)}
\eeq
and
\be\label{eq4-3}
\max\left\{\frac{d'(f(x), f(y))}{d_{G'}(x)}, \frac{d'(f(x), f(z))}{d_{G'}(x)}\right\}\leq c_1\max\left\{\frac{d(x, y)}{d_G(x)},\frac{d(x, z)}{d_G(x)}\right\}\leq c_1\vartheta_1\leq \frac{1}{2}.
\ee
It follows that
\beq\label{eq4-4}
\frac{1}{2}d_{G'}(f(x))\leq  d_{G'}(f(y))\leq 2d_{G'}(f(x)).
\eeq
By \eqref{eq4-3} and \eqref{eq4-4}, we get
\beqq\label{eq4-5}
d'(f(z), f(y))&\leq& d'(f(x), f(z))+d'(f(x), f(y))\\
&\leq& 2c_1\vartheta_1d_{G'}(f(x))
\leq 4c_1\vartheta_1d_{G'}(f(y))
< t_0d_{G'}(f(y)).
\eeqq

On the other hand, since $f^{-1}$ is $(\theta, t_0)$-relative, we have
$$\frac{d(y, z)}{d_G(y)}\leq c_1\frac{d'(f(y), f(z))}{d_{G'}(f(y))}.$$
Hence we obtain from \eqref{eq4-1} and \eqref{eq4-4} that
\beq\label{eq4-6}
\frac{d(y, z)}{d_{G}(x)}\leq 4c_1\frac{d'(f(y), f(z))}{d_{G'}(f(x))}.
\eeq
By \eqref{eq4-1}, \eqref{eq4-2a}, and \eqref{eq4-4}, we have
\beqq\label{eq4-7}
\frac{d'(f(y), f(z))}{d_{G'}(f(x))}\leq 4c_1\frac{d(y, z)}{d_G(x)},
\eeqq
which, together with \eqref{eq4-6}, implies
$$\frac{1}{4c_1} C_x d(y, z)\leq d'(f(y), f(z))\leq 4c_1  C_x d(y, z),$$
where $C_x=d_{G'}(f(x))/d_G(x)$.  The proof is complete.
\epf

\begin{lem}\label{lem3-5}
Let $0<\vartheta_1<1$ and $L_1\geq 1$. If $f$ is $(L_1, \vartheta_1)$-locally biLipschitz, then $f$ is $(\eta, q)$-locally quasisymmetric with $\eta(t)=c_2t$, where $c_2$ and $q$ depend only on $\vartheta_1$ and $L_1$.
\end{lem}
\bpf
Let $q=\vartheta_1$ and $\eta(t)=L_1^2 t$. Fix $u\in G$ and   three distinct points $x, y, z\in B(u, q d_G(u))$. By the assumption, we have
$$\frac{C_u}{L_1}d(x, y)\leq d'(f(x), f(y))\leq L_1C_u d(x, y),$$
and
$$\frac{C_u}{L_1}d(x, z)\leq d'(f(x), f(z))\leq L_1C_u d(x, z).$$
Hence we obtain
$$\frac{d'(f(x), f(y))}{d'(f(x), f(z))}\leq L_1^2 \frac{d(x, y)}{d(x, z)},$$
which completes the proof.
\epf

\begin{lem}\label{lem3-6}
Let $0<q<1$ and $c_2\geq 1$. If both $f$ and $f^{-1}$ are $(\eta, q)$-locally quasisymmetric with $\eta(t)=c_2t$, then $f$ is $\partial$-biLipschitz with data $L$ and $\lambda$, where $L$ and $\lambda$ depend only on $q$, $c_2$ and $c$.
\end{lem}
\bpf Let $q_1=\min\{1/(2+c), q/2\}$,
$L=8cc_2/q_1$ and $\lambda=q_1/(2cc_2)$. We need to prove that $f$ is $\partial$-Lipschitz with data $(L, \lambda)$, because the other direction follows from a similar argument. For all $x\in G$ and $y, z\in B(x, \lambda d_G(x))$, we have
$ 2^{-1}d_G(x)\leq d_G(y)\leq 2d_G(x).$
Let $u'\in \partial G'$ be such that $d'(f(x), u')\leq 2d_{G'}(f(x)).$ Since $X'$ is
$c$-quasiconvex, there exists a curve $\gamma_1'$ in $X'$ joining $u'$ and $f(x)$ such that
\beq\label{eq4-8b}
\ell(\gamma_1')\leq cd'(f(x), u')\leq 2cd_{G'}(f(x)).
\eeq
It follows from Lemma~C that
\beq\label{eq4-9}
B(x, q_1d_G(x))\subset G.
\eeq
So there is a point $u_1'$ belonging to $\gamma_1'\cap f(S(x, q_1d_G(x)))$. We see from \eqref{eq4-8b} that
\beq\label{eq4-10}
d'(f(x), u_1')\leq \ell(\gamma_1')\leq cd'(f(x), u')\leq 2cd_{G'}(f(x)).
\eeq
Because $f$ is $(\eta, q)$-locally quasisymmetric, we have
$$\frac{d'(f(x), f(y))}{d'(f(x), u_1')}\leq c_2\frac{d(x, y)}{d(x, u_1)}\leq c_2\frac{\lambda d_G(x)}{q_1d_G(x)}= \frac{c_2\lambda}{q_1},$$
where $u_1=f^{-1}(u_1')$. Using \eqref{eq4-10}, we obtain
$$d'(f(x), f(y))\leq \frac{c_2\lambda}{q_1}d'(f(x), u_1')\leq \frac{2cc_2\lambda}{q_1}d_{G'}(f(x))\leq \frac{1}{2}d_{G'}(f(x)),$$
which implies that
\beq\label{eq4-11}
\frac{1}{2}d_{G'}(f(x))\leq d_{G'}(f(y))\leq 2d_{G'}(f(x)).
\eeq

Choose $u_2'\in \partial G'$ such that
$d'(f(y), u_2')\leq 2d_{G'}(f(y)).$
The quasiconvexity of $X'$ ensures that there exists a curve $\gamma_2'$ in $X'$ joining $f(y)$ and $u_2'$ such that
\beq\label{eq4-12}
\ell(\gamma_2')\leq cd'(f(y), u_2')\leq 2cd_{G'}(f(y)).
\eeq
By \eqref{eq4-9}, there is a point $u_3'$ lying on $\gamma_2'\cap f(S(x, q_1d_G(x)))$. Then it follows from \eqref{eq4-12} that
$d'(f(y), u_3')\leq 2cd_{G'}(f(y)).$
Furthermore, by the local quasisymmetry of $f$ and by \eqref{eq4-11}, we have
\beqq
\frac{d'(f(y), f(z))}{d_{G'}(f(x))}&\leq& \frac{2d'(f(y), f(z))}{d_{G'}(f(y))}
\\ &\leq& 4c\frac{d'(f(y), f(z))}{d'(f(y), u_3')}\leq 4cc_2\frac{d(y, z)}{d(u_3, y)}
\\ &\leq& 4cc_2\frac{d(y, z)}{(q_1-\lambda)d_G(x)}\leq \frac{8cc_2}{q_1}\frac{d(y, z)}{d_G(x)}.
\eeqq
Hence we complete the proof of this lemma.
\epf

\noindent
\textbf{Proof of Theorem \ref{thm-2}.} The implications (\ref{en2-1})$\Rightarrow$(\ref{en2-2})$\Rightarrow$(\ref{en2-3})$\Rightarrow$(\ref{en2-1}) follow from Theorem \ref{thm-1} directly. From Lemma \ref{lem3-4} we obtain (\ref{en2-3})$\Rightarrow$(\ref{en2-4}). The implication (\ref{en2-4})$\Rightarrow$(\ref{en2-5}) follows from Lemma \ref{lem3-5}. By Lemma \ref{lem3-6}, we get the implication (\ref{en2-5})$\Rightarrow$(\ref{en2-1}).  The proof of Theorem \ref{thm-2} is complete.
\qed
\section{Proof of Theorem \ref{thm-3}}\label{sec-4}
\subsection{Gromov hyperbolic spaces}
Let $(X,d)$ be a metric space.
Fix a base point $w\in X$. For $x, y\in X$, the number
$$(x|y)_w=\frac{1}{2}\big(d(x, w)+d(y, w)-d(x, y)\big)$$
is called the Gromov product of $x, y$ with respect to $w$. We say that $X$ is a Gromov hyperbolic space (or a $\delta$-hyperbolic space) if there is a constant $\delta\geq 0$ such that
$(x|y)_w\geq \min\{(x|z)_w, (z|y)_w\}-\delta$
for all $x, y, z, w\in X$. For more information about Gromov hyperbolic spaces see \cite{BHK, Bu, ZP}. In particular, a direct computation gives the following elementary fact: For all $x,y,z,w,o\in X$, we have
\be\label{z-4}(x|y)_o+(z|u)_o-(x|z)_o-(y|u)_o
=(x|y)_w+(z|u)_w -(x|z)_w-(y|u)_w.\ee

\begin{defn}
Let $X$ be a proper (i.e., all closed balls are compact), geodesic, $\delta$-hyperbolic space, and $K\geq 0$. Following \cite{BHK}, we say that $X$ is $K$-roughly starlike with respect to a distinguished point $w\in X$ if for each $x\in X$, there is a geodesic ray $\gamma$ emanating from $w$ such that ${\rm dist}(x, \gamma)\leq K$.
\end{defn}

\vspace{8pt}
\noindent
{\bf Lemma E.} $($\cite{BHK}$)$ 
{\it
Let $(\Omega,d)$ be a locally compact and $A$-uniform metric space and $k$ its quasihyperbolic metric. Then there are constants $\delta,K\geq 0$ depending only on $A$ so that $(\Omega,k)$ is a proper, geodesic $\delta$-hyperbolic metric space. If, in addition, $\Omega$ is bounded, then $(\Omega,k)$ is $K$-roughly starlike with respect to $w$, where $w$ is a point in $\Omega$ satisfying $d_{\Omega}(w)=\max_{x\in \Omega}d_{\Omega}(x)$.
}

\vspace{8pt}

\subsection{Bonk-Heinonen-Koskela's uniformization}
Let $(\Omega, d)$ be a locally compact, rectifiably connected and incomplete metric space, where the identity mapping $\id: (\Omega, d)\to (\Omega, \ell)$ is a homeomorphism.  Let $k$ be the quasihyperbolic metric of $\Omega$ and $w\in \Omega$. In \cite{BHK}, Bonk et al. introduced a family of conformal deformations of $(\Omega,k)$ by the densities
$ \rho_{w,\varepsilon}(x)=e^{-\varepsilon k(x,w)}$ $(\varepsilon>0).$
For all $x$, $y\in \Omega$, we define
 \be\label{a-5} d_{w,\varepsilon}(x,y)=\inf\int_{\gamma} \rho_{w,\varepsilon} \; ds_k,\ee
where $ds_k$ denotes the arc-length element with respect to the metric $k$ and the infimum is taken over all rectifiable curves $\gamma$ in $\Omega$ with endpoints $x$ and $y$. Then $d_{w,\varepsilon}$  are metrics on $\Omega$, and we denote the resulting metric spaces by $\Omega_{w,\varepsilon}=(\Omega,d_{w,\varepsilon})$.

Now, we recall some useful facts concerning Bonk-Heinonen-Koskela's uniformization.

\vspace{8pt}
\noindent
{\bf Lemma F.} $($\cite{BHK}$)$
{\it If $(\Omega,k)$ is $\delta$-hyperbolic and $K$-roughly starlike with respect to $w\in\Omega$, then there are constants $A,C\geq 1,\varepsilon\in(0,1)$ that depend only on $\delta$ and a constant $M=M(\delta,K)\geq 1$ such that
\begin{enumerate}
\item[{\rm(a)}]\label{z-3a}   $\Omega_{w,\varepsilon}$ is $A$-uniform and bounded (with  diameter at most $2/\varepsilon$).

\item[{\rm(b)}]\label{z-3b}  The identity mapping $(\Omega,k)\to (\Omega,k_{w,\varepsilon})$ is $M$-biLipschitz, where $k_{w,\varepsilon}$ is the quasihyperbolic metric of $\Omega_{w,\varepsilon}$.

\item[{\rm(c)}]\label{z-3c} For all $x,y\in \Omega$, we have
 $$C^{-1}d_{w,\varepsilon}(x,y) \leq \varepsilon^{-1}e^{-\varepsilon (x|y)_w} \min\{ 1,\varepsilon k(x,y) \} \leq Cd_{w,\varepsilon}(x,y).$$
 \end{enumerate}
}

\vspace{8pt}

\begin{lem}\label{cl5-4}
Let $(\Omega, d)$ be a locally compact and $A$-uniform space. Suppose that $d_{w_0, \varepsilon}$ and $d_{w_1, \varepsilon}$ are metrics defined as in \eqref{a-5} with respect to the base points $w_0$ and $w_1$, respectively, where $\varepsilon=\varepsilon(A)\in (0, 1)$ is a parameter. Then the identity mapping
$(\Omega, d_{w_0, \varepsilon})\to (\Omega, d_{w_1, \varepsilon})$
is $\theta_1$-quasim\"obius, where $\theta_1$ depends only on $A$.
\end{lem}
\bpf It follows from Lemma~E, 
Lemma~F(c) 
and (\ref{z-4}) that there is a constant $C$ depending only on $A$ such that for any $x, y, z, u\in \Omega$,
\beqq
\frac{d_{w_1, \varepsilon}(x, y)d_{w_1, \varepsilon}(z, u)}{d_{w_1, \varepsilon}(x, z)d_{w_1, \varepsilon}(y, u)}
&\leq&
 C^4\frac{e^{-\varepsilon(x|y)_{w_1}}\min\{1, \varepsilon k(x, y)\}e^{-\varepsilon(z|u)_{w_1}}\min\{1, \varepsilon k(z, u)\}}{e^{-\varepsilon(x|z)_{w_1}}\min\{1, \varepsilon k(x, z)\}e^{-\varepsilon(y|u)_{w_1}}\min\{1, \varepsilon k(y, u)\}}\\
&=&
C^4\frac{e^{-\varepsilon(x|y)_{w_0}}\min\{1, \varepsilon k(x, y)\}e^{-\varepsilon(z|u)_{w_0}}\min\{1, \varepsilon k(z, u)\}}{e^{-\varepsilon(x|z)_{w_0}}\min\{1, \varepsilon k(x, z)\}e^{-\varepsilon(y|u)_{w_0}}\min\{1, \varepsilon k(y, u)\}}\\
&\leq&
C^8 \frac{d_{w_0, \varepsilon}(x, y)d_{w_0, \varepsilon}(z, u)}{d_{w_0, \varepsilon}(x, z)d_{w_0, \varepsilon}(y, u)},
\eeqq
which implies that the identity mapping $(\Omega, d_{w_0,\varepsilon})\to (\Omega,d_{w_1, \varepsilon})$ is $\theta_1$-quasim\"obius, where $\theta_1(t)=C^8t$.
\epf

\subsection{Proof of Theorem \ref{thm-3}}
In this subsection, we assume that $A\geq 1$, $0<q<1$, and $\eta:\,[0, \infty)\to [0, \infty)$ is a homeomorphism. Let $\Omega$ and $\Omega'$ be  locally compact and $A$-uniform spaces, and let $k$ and $k'$ be the quasihyperbolic metrics of $\Omega$ and $\Omega'$, respectively. Suppose that $f:\Omega\to \Omega'$ is a homeomorphism such that both $f$  and its inverse $f^{-1}$ are $(\eta, q)$-locally quasisymmetric.

\begin{lem}\label{cl5-1}
 The mapping $f$ is quasi-isometric with respect to the quasihyperbolic metrics.
 \end{lem}
\bpf For all $x, y\in \Omega$, we prove that there exist constants $L_1$ and $C_1$ such that
$$L_1^{-1}k(x,y)-C_1\leq k'(f(x), f(y))\leq L_1k(x, y)+C_1.$$
It suffices to check the right hand side inequality, because the left hand side inequality follows from a similar argument.  Moreover, \cite[Theorem 2.2]{Vai-3} tells us that we only need to show that if $x, y\in \Omega$ with $k(x, y)\leq t_1$, then we have
\beq\label{eq5-2}
k'(f(x), f(y))\leq 4A^2\log 2,
\eeq
where
$$q_1=\frac{q}{2}\,\,\,\;\;\mbox{and}\;\;\,\,\,
t_1=\min\left \{\log\left (1+\eta^{-1}\Big(\frac{1}{4A}\Big)q_1\right ),\log\Big(1+\frac{q}{2}\Big)
\right\}.
$$

To this end, we fix two points $x, y\in \Omega$ with $k(x, y)\leq t_1$. By Lemma~B(\ref{HRWZ-1}), 
we have
\beq\label{eq5-4}
d(x, y)\leq (e^{t_1}-1)d_{\Omega}(x).
\eeq

Let $u'\in \partial \Omega'$ be such that $d'(f(x), u')\leq 2d_{\Omega'}(f(x)).$
Because $\Omega'$ is $A$-uniform and locally compact, it follows from Lemma~E 
that  $(\Omega', k')$ is a proper, geodesic and Gromov hyperbolic space. Applying the Arzela-Ascoli Theorem, we know that there exists a quasihyperbolic geodesic ray $\gamma'$ in $\Omega'$ joining $f(x)$ and $u'$ that is an $A_1$-uniform curve  by using \cite[Theorem 2.10]{BHK}, where $A_1$ depends only on $A$. For convenience, we assume that $\gamma'$ is  an $A$-uniform curve by taking $A$ to be larger. Hence we obtain
\beq\label{eq5-5}
\ell(\gamma')\leq A d'(f(x), u')\leq 2Ad_{\Omega'}(f(x)).
\eeq
 Then there is a point $z'$ belonging to $\gamma'\cap f(S(x, q_1d_{\Omega}(x)))$. Let $z=f^{-1}(z')$.
Thus, by the local quasisymmetry of $f$, we see from \eqref{eq5-4} that
\beq\label{eq5-6}
\frac{d'(f(x), f(y))}{d'(f(x), z')}\leq \eta\Big(\frac{d(x, y)}{d(x, z)}\Big)\leq \eta\Big(\frac{e^{t_1}-1}{q_1}\Big)\leq \frac{1}{4A},
\eeq
where the last inequality follows from the choice of $t_1$.
By \eqref{eq5-5}, we have
$$d'(f(x), z')\leq \ell(\gamma')\leq 2 A d_{\Omega'}(f(x)),$$
which, together with \eqref{eq5-6}, shows that
$d'(f(x), f(y))\leq d_{\Omega'}(f(x))/2.$
Then it follows from \cite[(2.16)]{BHK} that
$$k'(f(x), f(y))\leq 4A^2\log\Big(1+\frac{d'(f(x), f(y))}{\min\{d_{\Omega'}(f(x)), d_{\Omega'}(f(y))\}}\Big)\leq 4A^2\log 2,
$$
which implies \eqref{eq5-2}. The proof of Lemma \ref{cl5-1} is complete.
\epf

For proving  Theorem \ref{thm-3}, we consider four cases in terms of the boundedness of the spaces $\Omega$ and $\Omega'$.

\begin{lem}\label{lem5-1}
Suppose that both $\Omega$ and $\Omega'$ are bounded. Then  $f$ is $\theta$-quasim\"obius, where $\theta$ depends only on $A$, $q$, and $\eta$.
\end{lem}
\bpf
Let $w_0\in \Omega$ and $w_1'\in \Omega'$ be such that
\beqq\label{eq5-1b}
d_{\Omega}(w_0)=\max_{x\in \Omega}d_{\Omega}(x)\;\;\mbox{and}\;\;d_{\Omega'}(w_1')=\max_{x'\in \Omega'}d_{\Omega'}(x'),
\eeqq
respectively. We first check that there is a constant $C$ depending on $A$ such that
\beq\label{eq5-9}
{\rm diam}(\Omega)\leq Cd_{\Omega}(w_0)\;\;\mbox{and}\;\;
{\rm diam}(\Omega')\leq Cd_{\Omega'}(w_1').
\eeq

For this, we only need to prove that
${\rm diam}(\Omega)\leq Cd_{\Omega}(w_0)$, because the second estimation follows from a similar argument. Take $u_1, u_2\in \Omega$ satisfying  $d(u_1, u_2)\geq {\rm diam}(\Omega)/2.$ Since $\Omega$ is $A$-uniform, we see that there is  an $A$-uniform curve $\gamma$ joining $u_1$ and $u_2$. Let $u\in \gamma$ be such that $\ell(\gamma[u_1, u])=\ell(\gamma[u, u_2])$. Therefore, we have
$\ell(\gamma)\leq 2Ad_{\Omega}(u)\leq 2Ad_{\Omega}(w_0),$
which guarantees that
${\rm diam}(\Omega)\leq 2\ell(\gamma)\leq 4Ad_{\Omega}(w_0).$
Hence, \eqref{eq5-9} holds.

It follows from Lemma~E 
that $(\Omega, k)$ and $(\Omega', k')$ are $\delta$-hyperbolic spaces, where $\delta$ depends only on $A$. Let $\Omega_{w_0, \varepsilon}=(\Omega, d_{w_0, \varepsilon})$ be the conformal deformation of $(\Omega, k)$ based at $w_0$ with parameter $0<\varepsilon\leq \varepsilon_0(A)$ which is defined as in \eqref{a-5}. Put $w_1=f^{-1}(w_1')$. Similarly, let $\Omega_{w_1, \varepsilon}=(\Omega, d_{w_1, \varepsilon})$ be the conformal deformation of $(\Omega, k)$ based at $w_1$ with parameter $\varepsilon$. It follows from Lemma~F 
that $\Omega_{w_0, \varepsilon}$ and $\Omega_{w_1, \varepsilon}$ are $A_1$-uniform, where $A_1$ depends only on $A$. Moreover, by \cite[(4.6)]{BHK} and Lemma~F(a), 
we have
\beq\label{eq5-10-1}
d_{w_0, \varepsilon}(w_0)\geq \frac{1}{\varepsilon e}\rho_{w_0, \varepsilon}(w_0)=\frac{1}{\varepsilon e}\geq \frac{1}{2e}{\rm diam}(\Omega_{w_0,\varepsilon})
\eeq
and similarly,
\beq\label{eq5-10}
d_{w_1, \varepsilon}(w_1)\geq \frac{1}{2e}{\rm diam}(\Omega_{w_1,\varepsilon}).
\eeq

To prove that $f:(\Omega,d)\to(\Omega',d')$ is $\theta$-quasim\"obius, we need the following:
\begin{enumerate}
\item[\textbf{I:}]  The identity mapping $\varphi_1:(\Omega, d)\to (\Omega, d_{w_0, \varepsilon})$ is $\eta_1$-quasisymmetric.

\item[\textbf{II:}] It follows from Lemma \ref{cl5-4} that the identity mapping $\varphi_2:(\Omega, d_{w_0, \varepsilon})\to (\Omega, d_{w_1, \varepsilon})$ is $\theta_0$-quasim\"obius with $\theta_0$ depending only on $A$.

\item[\textbf{III:}]  The induced mapping
$g:=f\circ \varphi_1^{-1}\circ \varphi_2^{-1}: (\Omega, d_{w_1, \varepsilon})\to (\Omega',d')$
is $\eta_2$-quasisymmetric.
\end{enumerate}
After that, we see that $f=g\circ \varphi_2\circ \varphi_1$ is $\theta$-quasim\"obius, because the composition of quasisymmetric mapping and quasim\"obius mapping is also quasim\"obius.

So in the remaining part, we only need to prove the assertions \textbf{I} and \textbf{III}.

For \textbf{I}, we check that the identity mapping $\varphi_1:(\Omega, d)\to (\Omega, d_{w_0, \varepsilon})$ satisfies the assumptions of Lemma~D. 
By (\ref{eq5-9}) and (\ref{eq5-10-1}), it suffices to demonstrate that both $\varphi_1$ and $\varphi_1^{-1}$ are locally quasisymmetric. Note that this desired claim follows from  
Lemma~F 
and Theorem \ref{thm-2}. Therefore, by  Lemma~D, 
the assertion \textbf{I} is valid.

Finally, it remains to prove \textbf{III}. We shall use Lemma~D 
to show that $g$ is $\eta_2$-quasisymmetric. From (\ref{eq5-9}) and (\ref{eq5-10}), we only need to verify that both $g$ and $g^{-1}$ are locally quasisymmetric. Since $f$  and $f^{-1}$ are both $(\eta, q)$-locally quasisymmetric, we are left to demonstrate that the identity mapping
$$\varphi_2\circ\varphi_1: (\Omega, d)\to (\Omega, d_{w_1, \varepsilon})$$
and its inverse are both locally quasisymmetric. This is because the composition of  locally quasisymmetric mappings  is also locally quasisymmetric, see \cite[Theorem 1.12]{HRWZ}.

This can be seen as follows. Because $\Omega'$ is a bounded $A$-uniform metric space, it follows from Lemma~E 
that $(\Omega', k')$ is $K'$-roughly starlike with respect to $w_1'$, where $K'$ depends only on $A$. Moreover, because the rough starlikeness is preserved under quasi-isometric mappings by \cite[Lemma 4.10]{ZP}, we see from Lemma \ref{cl5-1} that $(\Omega, k)$ is $K$-roughly starlike with respect to $w_1$, where $K$ depends only on $A$, $q$, and $\eta$. This fact, combined with Lemma~F(b), 
ensures that the identity mapping $\varphi_2\circ\varphi_1: (\Omega, d)\to (\Omega, d_{w_1, \varepsilon})$ is $M_1$-biLipschitz with respect to the quasihyperbolic metrics, where $M_1$ depends only on $A$, $q$, and $\eta$. It follows from Theorem \ref{thm-2} that both $\varphi_2\circ\varphi_1$ and its inverse are $(\eta_1, q_1)$-locally quasisymmetric, where $q_1$ and $\eta_1$ depend on $A$, $q$, and $\eta$.
The proof is complete.
\epf

Before proceeding, we introduce the sphericalization of metric spaces which were investigated by Bonk et al. in \cite{BK02} and Buckley et al. in \cite{BHX}.

\vspace{8pt}
\noindent
{\bf Lemma G.} 
{\it Suppose that $\Omega$ is an unbounded $A$-uniform space with $p\in \partial \Omega$, and that $(\Omega, \hat{d}_p)$ is the sphericalization of $(\Omega,d)$ associated to the point $p$.
Then the following statements hold:
\begin{enumerate}
\item[{\rm (1)}] $($\cite[Subsection 3.B]{BHX}$)$\label{en5-0} The space $(\Omega, \hat{d}_p)$ is a bounded metric space;
\item[{\rm (2)}] $($\cite[Subsection 3.B]{BHX} and \cite[Theorem 4.11]{BHX}$)$\label{en5-1}  The identity mapping $(\Omega, d)\to (\Omega, \hat{d}_p)$ is $\eta_0$-quasim\"obius and $M$-biLipschitz with respect to the quasihyperbolic metrics, where $\eta_0(t)=16t$ and $M=80A$;
\item[{\rm (3)}] $($\cite[Theorem 5.5]{BHX}$)$\label{en5-3}  The space $(\Omega, \hat{d}_p)$ is $A_1$-uniform, where $A_1$ depends only on $A$.
\end{enumerate}
}
\vspace{8pt}

\begin{lem}\label{lem5-2} Suppose that both $\Omega$ and $\Omega'$ are unbounded. Then $f$ is $\theta$-quasim\"obius, where $\theta$ depends only on $A$, $q$, and $\eta$.
\end{lem}
\bpf
Let $p\in \partial \Omega$ and $p'\in \partial \Omega'$. Suppose that $(\Omega,\hat{d}_p)$ and $(\Omega',\hat{d'}_{p'})$ are the sphericalizations of $(\Omega,d)$ and $(\Omega', d')$ associated to the points $p$ and $p'$, respectively. Lemma~G(1) 
ensures that $(\Omega, \hat{d}_p)$ and $(\Omega', \hat{d'}_{p'})$ are bounded metric spaces. Since $\Omega$ and $\Omega'$ are $A$-uniform, Lemma~G(3) 
implies that $(\Omega, \hat{d}_p)$ and $(\Omega', \hat{d'}_{p'})$ are $A_1$-uniform, where $A_1$ depends only on $A$.

By Lemma~G(3), 
we know that the identity mappings $\phi_1:\,(\Omega, d)\to (\Omega, \hat{d}_p)$ and $\phi_2:\,(\Omega', d')\to (\Omega', \hat{d'}_{p'})$  are $M$-biLipschitz with respect to the quasihyperbolic metrics, where $M$ depends only on $A$. Moreover,  it follows from Theorem \ref{thm-2} that $\phi_1$, $\phi_2$ and their inverses $\phi_1^{-1}$, $\phi_2^{-1}$ are $(\eta_1, q_1)$-locally quasisymmetric, where $\eta_1$ and $q_1$ depend only on $A$.
Thus, $f$ induces a mapping
$$g:=\phi_2\circ f\circ\phi_1^{-1}:\,(\Omega, \hat{d}_p)\to (\Omega', \hat{d'}_{p'}).$$
Since $f$  and $f^{-1}$ are both $(\eta, q)$-locally quasisymmetric, one finds that $g$ and its inverse $g^{-1}$ are $(\eta_2, q_2)$-locally quasisymmetric, where $\eta_2$ and $q_2$ depend only on $\eta$, $q$, and $A$. This is because the composition of  locally quasisymmetric mappings is also locally quasisymmetric, see \cite[Theorem 1.12]{HRWZ}.
Now, it follows from Lemma \ref{lem5-1} that $g$ is $\theta_1$-quasim\"obius, where $\theta_1$ depends only on $\eta$, $q$ and $A$.

Finally, Lemma~G(2) 
asserts that the identity mappings $\phi_1:\,(\Omega, d)\to (\Omega, \hat{d}_p)$ and $\phi_2:\,(\Omega', d')\to (\Omega', \hat{d'}_{p'})$  are $\theta$-quasim\"obius with $\eta_0(t)=16t$. Hence, $f$ is $\theta$-quasim\"obius with $\theta$ depending only on $\eta$, $q$, and $A$, because the composition of quasim\"obius mappings is also quasim\"obius. The proof of the lemma is complete.
\epf

\begin{lem}\label{lem5-3}
  Suppose that $\Omega$ is unbounded and $\Omega'$ is bounded. Then $f$ is $\theta$-quasim\"obius, where $\theta$ depends only on $A$, $q$, and $\eta$.
\end{lem}
\bpf
Let $p\in \partial \Omega$. Suppose that $(\Omega, \hat{d}_p)$ is the sphericalization of $(\Omega,d)$ associated to the point $p$. Since $\Omega$ is $A$-uniform, 
one observes from Lemma~G(1) 
and (3) that $(\Omega, \hat{d}_p)$ is a bounded $A_1$-uniform space, where $A_1$ depends only on $A$.

By Lemma~G(2), 
we know that the identity mapping $\phi_1:\,(\Omega, d)\to (\Omega, \hat{d}_p)$ is $M$-biLipschitz with respect to the quasihyperbolic metrics, where $M$ depends only on $A$.  Then it follows from Theorem \ref{thm-2} that $\phi_1$ and its inverse $\phi_1^{-1}$ are $(\eta_1, q_1)$-locally quasisymmetric, where $\eta_1$ and $q_1$ depend only on $A$. Therefore, $f$ induces a mapping
$$g:=f\circ\phi_1^{-1}:\,(\Omega, \hat{d}_p)\to (\Omega', d').$$
Since $f$  and $f^{-1}$ are both $(\eta, q)$-locally quasisymmetric, we see that $g$ and its inverse $g^{-1}$ are $(\eta_2, q_2)$-locally quasisymmetric, where $\eta_2$ and $q_2$ depend only on $\eta$, $q$, and $A$. This is because the composition of  locally quasisymmetric mappings is also locally quasisymmetric, see \cite[Theorem 1.12]{HRWZ}.
Now, it follows from Lemma \ref{lem5-1} that $g$ is $\theta_1$-quasim\"obius, where $\theta_1$ depends only on $\eta$, $q$ and $A$.

Finally, Lemma~G(2) 
tells us that the identity mapping $\phi_1:\,(\Omega, d)\to (\Omega, \hat{d}_p)$ is $\theta$-quasim\"obius with $\eta_0(t)=16t$. Hence, $f$ is $\theta$-quasim\"obius with $\theta$ 
depending only on $\eta$, $q$, and $A$, because the composition of quasim\"obius mappings is also quasim\"obius. The proof of the lemma is complete.
\epf

\begin{lem}\label{lem5-4}
Suppose that $\Omega$ is bounded and $\Omega'$ is unbounded. Then $f$ is $\theta$-quasim\"obius, where $\theta$ depends only on $A$, $q$, and $\eta$.
\end{lem}
\bpf
The proof of the lemma is similar to that of Lemma \ref{lem5-3} by exchanging the roles of $\Omega$ and $\Omega'$.
\epf

\noindent
\textbf{Proof of Theorem \ref{thm-3}.} Theorem \ref{thm-3} follows from  Lemmas \ref{lem5-1}\,--\,\ref{lem5-4}.
\qed

\bigskip
%


\end{document}